\documentclass[11pt,reqno]{amsart}
\usepackage{amsmath,amsthm}
\usepackage{amssymb,latexsym}
\usepackage{enumerate}

\usepackage[english]{babel}
\usepackage{booktabs}
\usepackage{url}

\usepackage[latin1]{inputenc}
\usepackage[T1]{fontenc}

\newcommand{\Z}{{\mathbb Z}}

\newcommand{\Q}{{\mathbb Q}}

\newcommand{\ta}{{\theta }}
\newcommand{\ent}{E_{n,\ta}}

\newcommand{\ttt}{E_{\ta}}

\theoremstyle{plain}
\numberwithin{equation}{section}
\newtheorem{thm}{Theorem}[section]
\newtheorem{theorem}[thm]{Theorem}
\newtheorem{lemma}[thm]{Lemma}

\title{On the high rank $\pi/3$ and $2\pi/3$-congruent number elliptic curves }

\author{A. S. Janfada}
\address{Department of Mathematics, Urmia University, Urmia, Iran}
\email{a.sjanfada@urmia.ac.ir; asjanfada@gmail.com}

\author{S. Salami}
\address{Instituto da Mathemçática  e Estatistica,
         UERJ, Rio de Janeiro,  Brazil}
\email{sajad.salami@ime.uerj.br}

\author[A. Dujella]{A. Dujella}
\address{
Department of Mathematics\\
University of Zagreb\\
Bijeni{\v c}ka cesta 30, 10000 Zagreb, Croatia}
\email{duje@math.hr}

\author[J. C. Peral]{J. C. Peral}
\address{
Departamento de Matematicas\\
Universidad del Pais Vasco\\
Aptdo. 644, 48080 Bilbao, Spain}
\email{juancarlos.peral@ehu.es}
\date{}
\keywords{$\ta$-congruent number, elliptic curve, Mordell-Weil rank}
\subjclass{11G05}

\begin{document}

%%==============================================================================
\begin{abstract}
Consider the elliptic curves    given by
 \[
 \ent :\quad
  y^2=x^3+2s n x^2-(r^2-s^2) n^2 x
  \]
 where $0 < \ta< \pi$,  $\cos(\ta)=s/r$ is rational with  $0\leq |s| <r$ and $\gcd
(r,s)=1$. These elliptic curves are related to the
$\ta$-congruent number problem as a generalization of the
 congruent number problem. For fixed $\ta$ this family corresponds
 to the quadratic twist by $n$ of the curve $\ttt: \,\,  y^2=x^3+2s  x^2-(r^2-s^2) x.$
We study two special cases
$\ta=\pi/3$ and $\ta=2\pi/3$. We have found a subfamily of $n=n(w)$ having rank at least
$3$ over $\Q(w)$ and a subfamily with rank $4$ parametrized by points of
an elliptic curve with positive rank.
We also found examples of  $n$ such that $E_{n, \ta}$
has  rank up to $7$ over $\Q$ in both cases. 
\end{abstract}

\maketitle

%%==========================================================================================
\section{Introduction}
\label{intro} The construction of  high rank elliptic curves is an important
problem concerning  elliptic curves. Dujella
\cite{duje} collected a list of  high rank elliptic curves
with prescribed torsion groups. The largest known rank, found by
 Elkies \cite{elk} in 2006, is $28$. In this work we  search for high ranks in the family of elliptic curves related with $\pi/3$ and $2\pi/3$ congruent problem.

 %----------------------------

 Let us  briefly describe the problem.
Consider  $0<\ta<\pi$  such that
$\cos(\ta)=s/r$ with $r$ and $s$ in $\Q$, $0\leq |s| <r$ and
$\gcd (r,s)=1$. A positive integer $n$ is called a
$\ta$-congruent number if there exists a triangle with rational
sides and area equal to $n{\alpha_\ta}$, where
${\alpha_\ta}=\sqrt{r^2-s^2}$. It is clear that if a positive integer $n$ is $\theta$-congruent,
then so is $n t^2$, for any  integer $t$, so we concentrate on square-free positive integers.

The problem of determining $\ta$-congruent numbers is
related to the problem of finding non-2-torsion points on the
family of elliptic curves which are called  $\ta$-congruent number elliptic curves,
\[
\ent: y^2=x^3+2snx^2-(r^2-s^2)n^2x,
\]
where $r$ and $s$ are as above, see \cite{topyui}. Observe that this curve is the quadratic twist by $n$ of the curve $E_{1, \ta}$.

This family of elliptic curves was introduced by Koblitz in \cite[Section I.2, Exercise 3]{kob1}, 
and systematically studied by Fujiwara \cite{fujw1, fujw2}. 
Let $E_{n,\ta}(\Q)$ be the group of rational points on $\ent$ and denote by $r_{\ta}(n)$
its (algebraic)  rank.

An ordinary congruent number is nothing but
a $\pi/2$-congruent number and hence a  congruent number elliptic curve is just a  $\pi/2$-congruent number elliptic curve. Rogers \cite{rgs1,rgs2} and
Dujella, Janfada and Salami
\cite{duj-js}, exhibited recently a  list of  congruent number elliptic curves
with $r_{\pi/2}(n)$ up to $7$.

We restrict our search   for high rank $\ta$-congruent number elliptic curves  to the  cases $\ta=\pi/3$ and $2\pi/3$.

In this paper we present a family of values of $n=n(w)$ such that the curves $E_{n(w),2 \pi/3}$   have  rank at least $3$ over $\Q(w)$. An equivalent result is  valid for the $\pi/3$ case.
We also exhibit examples of curves with rank up to $7$ in both cases, $\pi/3$ and $ 2\pi /3$.

 Yoshida
\cite{ShYsh1} proved  that $r_{\pi/3}(6)=1$,  $r_{\pi/3}(39)=2$ and also $r_{2\pi/3}(5)=1$,   $r_{2\pi/3}(14)=2$. 
These are the smallest positive integers corresponding to the given Mordell-Weil ranks. 
In this paper, we find the smallest
positive integers $n$ for which $r_{\pi/3}(n)=3,4,5$ in one case
and $r_{2\pi/3}(n)=3,4$ on the other (the result for $r_{\pi/3}(n)=4$ is conditional,
assuming the BSD and GRH).

 In our computations we use the Pari/Gp
software  \cite{pari}, William Stein's  SAGE software \cite{sage}
and Cremona's mwrank program \cite{crem1} and the program package Magma \cite{magma}.

%====================================================================
\section{Preliminary results}
In this section we recall some  results about $\ta$-congruent number elliptic curves,
in particular, a criterion
for a square-free positive integer to be a $\ta$-congruent number. This,  jointly with the subfamilies  mentioned before, are the starting point for our search of good candidates for high rank curves.

We use  the
Mestre-Nagao sum, the  Mestre's  conditional upper bound for the rank of
elliptic curves over $\Q$ and the root number  as sieving tools in order to reduce the size of the lists and selecting only the best candidates for high rank.  We briefly describe  these items  below.
%----------------------

It is known that for the usual congruent numbers there exist a close  relation with  elliptic curves, and
in fact the following classical result holds: $n$ is a congruent
number if and only if $r_{\pi/2}(n) >0$,  see e.g. \cite[Section I.9, Proposition 18]{kob1}.
 A similar theorem was proved by Fujiwara, see  \cite{fujw1}, for
$\ta$-congruent numbers.
\begin {theorem}\label{fuj}
Let $n$ be arbitrary   square-free positive integer and consider the
elliptic curve $E_{n,\ta}$ as above. Then
\begin{description}
\item[i)]  $n$ is a $\ta$-congruent number if and only if there exists a non-$2$-torsion point in  $\ent(\Q)$;
\item[ii)] for  $n\neq 1,2,3,6$, $n$ is a $\ta$-congruent number if and only if $r_{\ta}(n) >0$.
\end{description}
\end {theorem}

%-----------------------------
 Kan \cite{kan1} proved the following result which gives a family of  $\theta$-congruent
numbers for every $0< \ta< \pi$.
%------------------
\begin{lemma}\label{lem2}
A square-free positive integer $n$ is a $\ta$-congruent number if
and only if  $n$ is the square-free part of
\begin {equation}
\label{eq1} pq(p+q)(2rq+p(r-s)),
\end{equation}
for some positive integers $p$, $q$ with ${\rm gcd}(p,q)=1$.
\end{lemma}
%=================================================================

 Yoshida \cite{ShYsh1,ShYsh2} proved important results concerning  $\ta$-congruent numbers.  In particular, in \cite{ShYsh1} he gave the root numbers for the cases $\pi/3$ and $2\pi/3$
 (see Table \ref{t1:table}).
%begin{conjecture}  \mbox{ }
%\begin{itemize}
%\item [1)]$n$ is $\pi/3$-congruent number if $n\equiv
%6,10,11,13,17,18,21,22,23\ ({\rm mod}\ 24)$;
%\item [2)]$n$ is $2\pi/3$-congruent number if $n\equiv
%5,9,10,15,17,19,21,22,23\ ({\rm mod}\ 24)$.
%\end{itemize}
%\end{conjecture}

\begin{table}[htbp!]
\centering
\caption{Root-numbers}\label{t1:table}
\begin{tabular}{@{} l c c c @{}}
\toprule
       &       &$2 \pi/3$ & $\pi/3$\\
\midrule
$n\equiv 1, 2, 3, 6, 7,11,13,14,18$&   (mod $24$)     & $+ 1$   & $-1$      \\
$n\equiv 5,9, 10,15, 17, 19, 21, 22, 23$ &(mod $24$)  & $-1$    & $+ 1$    \\
\bottomrule
\end{tabular}
\end{table}
%=================================================================

Now we recall the Mestre-Nagao sum for an elliptic curve $E$ over $\Q$. Reduce $E$ modulo
a prime $p$ and suppose that  $N_p$  is the number of points on $E$ with coordinates on $F_p$.
 For any  positive integer $t$,
 let ${\bf P}_t$ be the set of all primes less than $t$ and $a_p=p+1 -N_p$.
 The  Mestre-Nagao sum is defined by
\[
S(t,E)=\sum_{p\in {\bf P}_t}^{}\left(1- \frac{p-1}{N_{p}}\right)\log p
=\sum_{p\in {\bf P}_t}^{} \frac{-a_p+ 2}{ N_{p}}\log p.
\]
It is experimentally known \cite{Me1,nag2}  that high
rank curves have large values $S(t,E)$. We cite \cite{camp} for a
heuristic argument which links the Mestre-Nagao sum to
the  Birch and Swinnerton-Dyer conjecture \cite{BSD}.

%------------------------

Now we    describe the   Mestre's conditional
(assuming the Birch and Swinnerton-Dyer conjecture and GRH)
upper bound (see \cite{Me2,djb}) for the rank of an elliptic curve over $\Q$. Let $E$ be an elliptic curve with conductor  $N$. For an integer $m\geq 1$, let
 \begin{displaymath}
  b(p^m)=\left\{ \begin{array}{ll}
   0 & \textrm{if $p \vert N$},\\
   \alpha_p^m + {\alpha'}_p^m & \textrm{if $ p\not\hspace{-.02cm} \vert N$},
 \end{array} \right.
 \end{displaymath}
where $\alpha_p$ and $\alpha'_p$ are the roots of $x^2-a_px+p$. Let
 \begin{displaymath}
  F(x)=\left\{ \begin{array}{ll}
   (1-x)\cos(\pi x)+\sin(\pi x)/\pi & \textrm{if $x \in [0,1]$},\\
   0 & \textrm{if $x>1$}.
 \end{array} \right.
 \end{displaymath}
Take  a positive real
number $\lambda$  and write
  \[
 M(\lambda)=2\big( \log(2\pi)+\int_0^{\infty}\left(F(x/\lambda)/(e^x -1)- e^{-x}/x \right)dx \big).
 \]
The  Mestre's conditional upper bound   for
the rank of $E$  is defined as
\[
M(\lambda,E)=\frac{\pi^2}{8\lambda} \Big ( \log(N) - 2 \sum_{p^m\leq e^\lambda} b(p^m) F(m\log(p)/\lambda)
\frac{\log(p)}{p^m} - M(\lambda)  \Big ).
\]

%=================================================================

%%%======================

\section{A family with generic rank at least $3$}
\label{fam}
\subsection{Twists}
Observe that, once $\ta$ is fixed, the curve
 \[
 \ent:\quad y^2=x^3+2s n x^2-(r^2-s^2) n^2 x
 \]
 is the quadratic twist with parameter $n$ of the curve $E_{1, \ta}: \quad  y^2=x^3+2s  x^2-(r^2-s^2) x.$ General results  about twists can be applied for any $\ta$ and  we can find families of rank at least $2$ over $\Q(r,s)$ by direct applications of results given in Mestre \cite{M} or Rubin and Silverberg \cite{RS1}, \cite{RS2} (see also \cite{GM,ST}).

 In our particular cases,  $\ta=\pi/3$ corresponds to $s=1$ and $r=2$ and
 $\ta=2\pi/3$ to $s=-1$ and $r=2$, we are lead to study the quadratic twists of the curves
\begin{align*}
 E_{\pi/3}:&\quad y^2=x^3+ 2 x^2- 3 x\\
 E_{ 2 \pi/3}:&\quad y^2=x^3- 2 x^2- 3 x
\end{align*}
Each curve is the twist  of the other by $-1$ so their twists can be studied jointly.

 \subsection{A family of twists for $\ta=2 \pi/3$ with rank $\geq 3$}
 \subsubsection{Rank  $1$}
 We start with the twists of the  curve
  \[
   E_{ 2 \pi/3}:\quad y^2=x^3- 2 x^2- 3 x
   \]
    with parameter $(u + a) (u + b) (u + c)$, so we have the family of twists $ y^2= x^3+ A x^2+ B x$ where
 \begin{align*}
 A & = - 2 (u + a) (u + b) (u + c),\\
 B & = - 3 (u + a)^2 (u + b)^2 (u + c)^2.
 \end{align*}
 Now we impose $-(b + u) (c + u)^2$ as the $x$-coordinate of a new point. This is the same as choosing
 \[
 c=\frac{-3 a - 4 u + a b w^2 + a u w^2}{1 + b w^2 + u w^2}
 \]
 With this choice we get a family of twists with rank at least $1$  over $\Q(b, u, w)$ which, after clearing denominators, can be written as $ y^2= x^3+ A_1 x^2+ B_1 x$  with
 \begin{align*}
 A_1 & = -2 (b + u) (-3 + b w^2 + u w^2) (1 + b w^2 + u w^2),\\
 B_1 & = -3 (b + u)^2 (-3 + b w^2 + u w^2)^2 (1 + b w^2 + u w^2)^2.
 \end{align*}
 The $x$-coordinate of the infinite order point is $x_1=-(b + u) (-3 + b w^2 + u w^2)^2  $.
  \subsubsection{Rank  $2$}
 We proceed by forcing $3 (b + u) (1 + b w^2 + u w^2)$ as the $x$-coordinate of a new point in the previous rank $1$ family of twists. For this purpose it is enough to choose
 \[
 b=-\frac{-4 + u^2 + u w^2 + u^3 w^2}{(1 + u^2) w^2}.
 \]
 Now the  new family of twist can be written as $ y^2= x^3+ A_2 x^2+ B_2 x$  with
 \begin{align*}
 A_2 & = -10 (-2 + u) (2 + u) (-1 + 2 u) (1 + 2 u) (1 + u^2),\\
 B_2 & = -75 (-2 + u)^2 (2 + u)^2 (-1 + 2 u)^2 (1 + 2 u)^2 (1 + u^2)^2.
 \end{align*}
 The $x$-coordinates of the two infinite order points are
  \begin{align*}
 x_1 & = (-2 + u) (2 + u) (-1 + 2 u)^2 (1 + 2 u)^2 (1 + u^2),\\
 x_2 & = -15 (-2 + u) (2 + u) (1 + u^2)^2.
 \end{align*}
 These two points are independent, so the new family has rank  at least $2$ over $\Q(u)$.
  \subsubsection{Rank  $3$} Finally we  choose
  \[
  u=-\frac{70 - 10 w + w^2}{3 (5 + w^2)}
  \]
  in order to get  $5 (-2 + u)^2 (-1 + 2 u)^2 (1 + u^2)$ as $x$-coordinate of a new point in the rank $2$ family. In this way we get $ y^2= x^3+ A_3 x^2+ B_3 x$  with
 \begin{align*}
 A_3 & = -2 (-5 + w) (-2 + w) (4 + w) (25 + w) (31 - 4 w + w^2)\\
     & \quad \quad \quad \quad \quad \ \  (100 - 10 w + 7 w^2) (1025 - 280 w + 66 w^2 - 4 w^3 + 2 w^4),\\
 B_3 & = -3 (-5 + w)^2 (-2 + w)^2 (4 + w)^2 (25 + w)^2 (31 - 4 w + w^2)^2 \\
	   & \quad \quad \quad \quad \ \  (100 - 10 w + 7 w^2)^2 (1025 - 280 w + 66 w^2 - 4 w^3 + 2 w^4)^2.
 \end{align*}
 The $x$-coordinates of the three independent  points are given by
 \begin{align*}
  x_1 & = \frac{-1}{9 (5 + w^2)^2} (-5 + w)^2 (-2 + w) (4 + w) (25 + w)^2 (31 - 4 w + w^2)^2 \\
	    & \quad \quad \quad \quad \quad \ \, (100 - 10 w + 7 w^2) (1025 - 280 w + 66 w^2 - 4 w^3 + 2 w^4),  \\
  x_2 & = 3 (-2 + w) (4 + w) (100 - 10 w + 7 w^2) \\ 
			&	\qquad \qquad \qquad \qquad \qquad \qquad \quad \, (1025 - 280 w + 66 w^2 - 4 w^3 + 2 w^4)^2,  \\ 
  x_3 & = (31 - 4 w + w^2)^2 (100 - 10 w + 7 w^2)^2 \\
	    & \qquad \qquad \qquad \qquad \qquad \qquad \quad \ \ (1025 - 280 w + 66 w^2 - 4 w^3 + 2 w^4).
 \end{align*}
 For  $w=10$, after reducing coefficients, we get the  rank $3$ curve given by $y^2= x^3-442 x^2-146523 x$.  The specialized points are   
 \begin{gather*}
 P_1=\{-2873/81, 1562912/729\}, \
 P_2=\{867,13872 \}, \\
 P_3=\{ 2873/4, 48841/8\}.
 \end{gather*}
 A calculation with mwrank \cite{crem1} shows that these three points are independent.  An argument of specialization \cite{S} proves that this family has rank at least $3$ over $\Q(w)$.

 Observe that the parameter for the rank $3$ family of twists can be made both positive  and negative for  infinitely  many values of $w$, so we get  a family of rank $3$ twist for both $2 \pi/3$ and $\pi/3$ congruent number problem.

\subsection{A subfamily with rank  $\geq 4$}
We can find a subfamily with rank $\geq 4$ in the family $ y^2= x^3+ A_3 x^2+ B_3r x$
by forcing
$$ -\frac{1}{4}(-5 + w)^2 (-2 + w)^2 (4 + w)^2 (25 + w)^2 (31-4 w+w^2)(100
- 10 w + 7 w^2) $$
to be the $x$-coordinate of a point on the curve.
We get the condition
\begin{equation} \label{quartic}
25 w^4-26 w^3+699 w^2-3770 w+13300 = z^2.
\end{equation}
It can be transformed to the elliptic curve
$$ Y^2 = X^3+X^2-17220X-352800 $$
with positive rank (rank is equal $2$ with generators $[255, 3450]$, $[-22, 126]$,
corresponding to the points $(w,z)=[315/74, 695275/5476]$, $[8, -342]$ on the quartic
(\ref{quartic})). Hence, we get infinitely many rational parameters $w$ for which
the curve $ y^2= x^3+ A_3 x^2+ B_3 x$ has the rank $\geq 4$.

%================================================================================
%=======================================================================================

\section{Strategies and  results}
\label{sec6}
\subsection{General setting}
\par Now we attempt to find high rank elliptic curves $E_{n,\ta}$  in two cases $ \ta=\pi/3 $ and
$2\pi/3$.  We will use the expression (\ref{eq1}) and  the families given in the previous section as sources for good candidates for high rank curves.
 We shall use the  following notations:
$r_\ta(n)$ for  the rank and
 $s_{\ta}(n)$ for the $2$-Selmer rank (see e.g. \cite{duj-js}),
 which is an upper bound for the rank; that is $r_\ta (n) \leq s_{\ta}(n)$.

 We proceed in three steps, depending on
the range and the form of the square-free positive integers $n$.

{\bf Step (I)} In this step we take all the square-free positive integers $n \leq 5\times 10^6$. By a direct computation with
mwrank, we find the 2-Selmer rank of $E_{n,\ta}$ for all square-free $n$ in that range and  in each case $\ta=\pi/3$ and $2\pi/3$.
Our computations show that there are no integers $n$ with $s_\ta(n)\geq 6.$
Table \ref{t:table2} presents the distribution of the number of these square-free integers according to  the  values of
$s_\ta(n)$.
Finally, we compute directly  rank $r_\ta(n)$ with mwrank to
  find the smallest $n$'s with
$ r_{\pi/3}(n)=3,4, 5$ as well as the smallest $n$'s with $ r_{2\pi/3}(n)=3, 4$.

\begin{table}[htbp!]
\centering
\caption{ Distribution of $s_\ta(n)$ }\label{t:table2}
\small
\begin{tabular}{@{} c c c c c c c c c @{}}
\toprule
$s_\ta(n)$ & $0$ &  $1$ & $2$ & $3$ &  $4$ & $5$ & $\geq6$   &  Total\\
\midrule
 $\ta=\pi/3$ & 783043 &  1401045 & 734290 & 116158 &  5045 & 52 &   0 & 3039633 \\
 \hline
 $\ta=2\pi/3$ & 760511 &  1374165 & 751192& 144641 &  9038 & 86 &  0 & 3039633 \\
\bottomrule
\end{tabular}
\end{table}

%=======================================================================================

 {\bf Step (II)}  We  consider all  square-free $\ta$-congruent numbers  $n> 5\times 10^6$ of the
form (\ref{eq1}) in Lemma \ref{lem2} with $1<p, q \leq 10^4$,
${\rm gcd}(p,q)=1$, and having at least $4$  odd prime factors. We
get a list  with more than $7\times 10^6$
elements for each of the cases $\ta=\pi/3$ and $\ta=2\pi/3$.
Using the Mestre-Nagao sum, we reduce by Pari/Gp program the length of this list. In fact, we choose
the $n$ with
 $S(10^3, E_{n,\ta})>15,  S(10^4, E_{n,\ta})>20, S(10^5, E_{n,\ta})>40,$
for which  $s_{\pi/3}(n) \geq 6$, and $s_{2\pi/3}(n) \geq 5$. After computing the
values of $r_\ta(n)$ for these candidates by mwrank, we finally select the
$n$ with $r_{\pi/3}(n)=6,7$ and the $n$ with
$r_{2\pi/3}(n)=5,6$. In the cases in which  mwrank do not give  exact value  $r_\ta(n)$  we compute the  Mestre's conditional upper bound $M(\lambda,E_{n,\theta})$ for $r_\ta(n)$ with
 $15\leq \lambda < 24$.

 {\bf Step (III)}   In this part we use the families in section \ref{fam}  in order to search for good candidates for high rank. Since curves in the families with rank $3$ and $4$ have large coefficients,
 we find the family with rank $2$ the most suitable for our purpose.
 The search for rank $6$ curves is conducted upon the rank $2$ family with  $u=p/q$ for $1 < p < q < 4000$  with sieving conditions $S (523, E_{n,\ta}) > 18$,  $S (1979, E_{n,\ta}) > 28$ and the Selmer rank  $\geq 6$.

 The search for rank $7$  is made in the same family of twists  with $u=p/q$ for $1 < p < q < 13000$ with the following conditions, root number equal to  $-1$, $S (523, E_{n,\ta}) > 20$,
 $S (1979, E_{n,\ta}) > 30$ and the Selmer rank  $\geq 7$.  The ranks are calculated with mwrank.  For the case $2\pi/3$ for $p=4127$ and  $q=10004$, i.e. for $n=12748697412909916241$ the corresponding curve has rank $7$. In this case, the direct application of mwrank gives  only $6\le$ rank $\le 7$, but
  applying mwrank to an isogenous curve give the seventh independent point.

\medskip

 %%%Here details................................

 In the next  subsections, we collect the results.
We  find the smallest integers  $n$ such that $r_{\pi/3}(n)=3,4,5$ and $r_{2\pi/3}(n)=3,4$,
and we exhibit examples of curves with rank up to $7$ in both cases.

%%In each sub case, using the LLL-reduction algorithm, we exhibit the independent set of rational points
%with the smallest height on the corresponding  $\ta$-congruent number elliptic curve $E_{n,\ta}$.
%Table 2 shows the final results.

%---------------------------------------------
\subsection{The case $\ta=\pi/3$}
\underline{Rank 3}: The integers $407$ and $646$ are the two
smallest ones among $116158$ integers $n$ less than $5\times
10^6$ with $s_{\pi/3}(n)=3$. We have
$r_{\pi/3}(646)=3$, while for
 $n=407$ Magma gives that the analytic rank is $1$, so by
  Kolyvagin's theorem $r_{\pi/3}(407)=1$. Therefore, the value $n=646$ is
 the minimum value producing  a curve with rank $3$.

%-------------------
 \underline{Rank 4}:
The smallest $n$ that we have found with rank $4$ is $n=172081$.
There are $63$ integers $n$ less than
$172081$  with $s_{\pi/3}(n)=4$.  For $29$ cases mwrank gives $0 \leq
r_{\pi/3}(n) \leq 4$, and for all these cases the $4$-descent implemented in Magma gives
that the rank is $\leq 2$.
In the remaining $34$ cases, mwrank gives $2 \leq
r_{\pi/3}(n) \leq 4$. In the most of these cases the $4$-descent shows that rank is equal to $2$.
However, in three cases: $n=31622$, $143222$, $150866$, we are not able to show that rank $<4$
unconditionally. In these cases, we use Mestre's conditional upper bound (with $\lambda=11$),
which gives $r_{\pi/3}(n) \leq 2$, so $r_{\pi/3}(n)=2$
(conditionally). Thus, the value $n=172081$ is, conditionally (assuming BSD and GRH),
the minimum value giving a curve with rank $4$.

%--------------------------
 \underline{Rank 5}: The direct computation shows that $n=221746$ is
the smallest among $52$ integers $n$ in the observed range
with $s_{\pi/3}(n)=5$, and since $r_{\pi/3}(221746)=5$, $n=221746$ is the smallest
positive integer giving rank $5$.

%--------------------------------
 \underline{Rank 6}: The smallest $n$ that we have found with rank $6$  is $n=11229594411$.
 We do not know if it is the smallest one with this property.
 The values of $n$ given in Table \ref{t:table3} also give curves with rank $6$.
 \begin{table}[htbp!]
\centering
\caption{Case $\pi/3$. Other $n$ with rank $6$}\label{t:table3}
\small
\begin{tabular}{@{} l l l l @{}}
\toprule
40004232681,&
158763281079,&
167514827545,& \\
198606002595,&
251819173095,&
271314827665,& \\
3302971161265,&
3492293850595,&
5144668978371,& \\
6634009064865,&
17073273800095,&
40582123000419,& \\
45563330326345,&
7658263493840940211& & \\
     \bottomrule
\end{tabular}
\end{table}

  %-------------------------------------
 \underline{Rank 7}:
The only $n$ that we have found giving rank $7$    is $n=365803464586$. We do not know if it is the smallest one.

%%%%%%%%%%%%%%%%%%%%%%%%%%%%%%%%%%%%%%%%%%%%%%%%%%%%%%%%%%%%%%%%%%%%%%%%%%%%%%%%%%%
\begin{table}[htbp!]
\label{tab2}
\caption{Ranks in the cases $\ta =\pi/3$}
\centering
\begin{tabular}{c}
{\tiny
\begin{tabular}{|c|c|l|}
\hline
  & & \\
 \scriptsize{ $r_{\ta}(n)$} &  \scriptsize{$n$} & \scriptsize{ Generators of $E_{n,\ta}: y^2=x^3+2snx^2-(r^2-s^2)n^2x$}  \\
                 &                 & \\
\hline
   & & \\
 3 & 646 & \verb [-722,34656] , \verb [6137,521645] , \verb [-1216,40432]  \\
   & & \\
\hline
   & & \\
 4 & 172081 & \verb [-505141,-61627202] ,  \verb [-58621,-78669382] , \\
   &        & \verb [-440076,-143244738] ,  \verb [224175,92987790]  \\
   & & \\
\hline
   & & \\
 5 & 221746 & \verb [345450,207822720] ,  \verb [-15792,49357896] , \verb [994896,1130036040] , \\
   &        & \verb [-13254,-45063600] ,  \verb [-386575,-255989965]   \\
   & & \\
\hline
   & & \\
 6 & 11229594411 & \verb [904103532759/25,-992069570757491352/125] , \\
   &             & \verb [1541731888897/16,2090318638263775025/64] , \\
   &             & \verb [265444083202036/2025,4636387440736982658134/91125] , \\
   &             & \verb [719501508201/64,40873417425022581/512] , \\
   &             & \verb [13006760076899764/269361,1693181585331404000267498/139798359] , \\
   &             & \verb [50286669020153449/278784,11896090671289659453790795/147197952]   \\
   & & \\
\hline
   & & \\
 7 & 365803464586 & \verb [433764757524,212456676940982628] , \\
   &              & \verb [1291274050073,-1689545579159165609] , \\
   &              & \verb [-59335333874904423/3644281,-570541659890431976790514695/6956932429] , \\
   &              & \verb [11954902524369/4,-45277466996084516865/8] , \\
   &              & \verb [2138828658027602/5329,56890395483549429623312/389017] , \\
   &              & \verb [786769181014433554/80089,721982407380536692088852160/22665187] ,  \\
   &              & \verb [-562236028164373765342/540237049,3617165210435366625559445197360/12556729729907]  \\
   & & \\
\hline
\end{tabular}} \\

\end{tabular}
\end{table}

%--------------------------------
\subsection{The case $\ta=2\pi/3$}
%--------------------------------
 \underline{Rank 3}: The smallest $n$ with $s_{2\pi/3}(n)=3$ is $n=221$.
 Since $r_{2\pi/3}(221)=3$, we conclude that $n=211$ is the smallest $n$ for which the rank is $3$.

%--------------------------------
\underline{Rank 4}:  The smallest $n$ that we have found with rank $4$ is $n=12710$.
There are two smaller positive integers with Selmer rank equal to $4$ ($n=4718$ and $n=6398$)
but having analytic rank $0$,
so by Kolyvagin's theorem the algebraic rank is also $0$. Thus the minimality of $n=12710$ follows.

%--------------------------------
 \underline{Rank 5}: The smallest $n$ that we have found with rank $5$ is  $n=16470069$.   We do not know if it is the smallest one with this property.

%--------------------------------
 \underline{Rank 6}: We have found several   positive integers  $n$
 with $r_{2\pi/3}(n)=6$ where $n=456249066$ is the
 smallest one.  Other  values are given in Table \ref{t:table4}.
 \begin{table}[htbp!]
\centering
\caption{ Case $2 \pi/3$. Other $n$ with rank $6$}\label{t:table4}
\small
\begin{tabular}{@{} l l l @{}}
\toprule
764046470,& 902472906,&  5062245006,\\ 
9667090290,& 11801899970,&19969987310,\\ 
20240772006,& 23819599518,& 24080567966,\\ 
30834423438,&39360775454,& 58181539130,\\ 
64256704710,& 98708770590,& 106366008126,\\ 
148280772990,& 181684390314,& 292826163630,\\ 
309000045354,& 333515184002,& 554883184814,\\ 
653918457570,& 685374515826,& 713465075246,\\ 
860842004286,& 1185986591790,&  1248260820170,\\ 
2004510092970,& 2743972777910,& 10745486363210,\\ 
55967962170246,& 90952836208430,& 104732378607110,\\ 
177348563238770,& 219163751391326,& 1459584795789354,\\ 
29410732919116094,& 40315634933149394,& 30375400815771401390 \\
\bottomrule
\end{tabular}
\end{table}

%%%%%%%%%%%%%%%%%%%%%%%%%%%%%%%%%%%%%%%%%%%%%%%%%%%%%%%%%%%%%%%%%%%%%%%%%%%%%%%%%%%
\begin{table}[htbp!]
\label{tab2a}
\caption{Ranks in the cases $\ta =2\pi/3$}
\centering

\begin{tabular}{c}
{\tiny
\begin{tabular}{|c|c|l|}
\hline
   & & \\
  \scriptsize{$r_{\ta}(n)$} &  \scriptsize{$n$} & \scriptsize{Generators of $E_{n,\ta}: y^2=x^3+2snx^2-(r^2-s^2)n^2x$}  \\
                 &                 & \\
\hline
   & & \\
 3 & 221 & \verb [-204,1734] , \verb [-169,2704] ,  \verb [4131,-249696]  \\
   & & \\
\hline
   & & \\
 4 & 12710 & \verb [-310,384400] ,  \verb [-9920,-1153200] , \\
   &       & \verb [48050,5381600] ,  \verb [76880,16337000]  \\
   & & \\
\hline
   & & \\
 5 & 16470069 & \verb [-3115959/4,-198146948769/8] , \\
   &          & \verb [-16255958103/1024,-813789518594283/32768] , \\
   &          & \verb [118172745075/1849,-21701053829180880/79507] , \\
   &          & \verb [174895662711/3481,-10850526914590440/205379] , \\
   &          & \verb [18013358979/361,-275820552686448/6859]  \\
   & & \\
\hline
   & & \\
 6 & 4562490669 & \verb [1372171206,2930957696016] , \verb [24303608784,3714988879700280] , \\
   &            & \verb [1677715326,-33259028622624] , \verb [3635049873,-183588193835865] , \\
   &            & \verb [27273656667348/18769,39342846732689875284/2571353] , \\
   &            & \verb [36967427406/25,2217080599939296/125]  \\
   & & \\
\hline
\end{tabular}} \\

\end{tabular}
%\caption{Rank records in two cases $\ta =\pi/3$ and $\ta =2\pi/3$  }
\end{table}
%+++++++++++++++++++++++++++++++++++++++++++++++++

 \underline{Rank 7}: The integer  $n=12748697412909916241$
 with $r_{2\pi/3}(n)=7$ has been found within the family of rank $2$ of section \ref{fam}. It correspond to  $u=\frac{4127}{10004}$ in such family. The data for this curve are too large to fit in the table, so we give them here. The rank and independent points were found by applying mwrank \cite{crem1}
 to one of its $2$-isogenous curves. The curve is
 \begin{multline*} 
y^2 = x^3-25497394825819832482x^2 \\
- 487587857177807974195124652448906710243x 
\end{multline*}
and $x$-coordinates of $7$ independent points are:
{\tiny
\begin{multline*}
  1. \quad -3478204633589378700, \\
  \shoveleft{2. \quad -11685945449719133341,} \\
  \shoveleft{3. \quad -6574179551855299730183742058990161459509481/575598836877796985970025,} \\
  \shoveleft{4. \quad 2582493196592574693159131199086103504610591321/64687220044469657223311844,} \\
  \shoveleft{5. \quad 937805074272703399240860666902959419125740930561/18861375626453019864153493504,} \\
  \shoveleft{6. \quad 11046597353655160109746607088518021300436780651...}\\
  \shoveleft{\quad \quad                        9686653884678929416873279865087505494460413809/} \\
        \shoveright{ 11748335750378251588082756719839642493430053195853237296682112610898176,} \\
  \shoveleft{7. \quad 604677170128352690791956712052363555568843154944505048349318...}\\
  \shoveleft{\quad \quad                        84615838950389619168774330021927698790345943423897889681/} \\
        \shoveright{ 2423301702231258764005536765918408361894206534...}\\
	      \shoveright{57595172659874108380495724928358099557362727572100}\\
\end{multline*} }%


\end{document}